\documentclass[12pt]{amsart}
\usepackage{amssymb,amscd}
\usepackage[bitstream-charter]{mathdesign}
\DeclareSymbolFont{cmsymbols}{OMS}{cmsy}{m}{n}
\SetSymbolFont{cmsymbols}{bold}{OMS}{cmsy}{b}{n}
\DeclareSymbolFontAlphabet{\mathcal}{cmsymbols}

\usepackage[%hypertex,
unicode=true,
plainpages = false,
pdfpagelabels,
bookmarks=true,
bookmarksnumbered=true,
bookmarksopen=true,
breaklinks=true,
backref=false,
colorlinks=true,
linkcolor = blue,
urlcolor  = blue,
citecolor = red,
anchorcolor = green,
hyperindex = true,
hyperfigures
]{hyperref}
\hypersetup{
	pdftitle={Linear algebra, theory and algorithms},
	pdfauthor={V. Mikaelian},
	pdfsubject={Linear Algebra}
}

\title[Higman operations on integer sequences]{\large Higman operations on integer sequences, and embeddings of recursive groups}

\author[V.\,H. Mikaelian]{V.\,H. Mikaelian\\ \lowercase{ v.mikaelian@gmail.com }}

\begin{document}

\maketitle

{\small\noindent This is an extended version of summary of the talk at the \textit{International Conference on Group Theory 
in honor of Victor Mazurov on the occasion of his 80th birthday}, 
Novosibirsk, 2--8 July 2023. The concise version of this report can be found in the talks and communications band at the Conference site:\\
\href{https://sites.google.com/view/mazurov2023}{\tt \footnotesize https://sites.google.com/view/mazurov2023}}

\bigskip
\noindent
The objective of the current talk is to present our very recent article \cite{The Higman operations and  embeddings} in which we for various types of recursive groups discuss possibilities of their \textit{explicit} embeddings into finitely presented groups.

\section{Higman's embedding theorem}
\noindent
In 1961 Higman proved that \textit{a finitely generated group can be embedded in a finitely presented group if and only if it is recursively presented} \cite{Higman Subgroups of fP groups}.
Higman's work is based on specific recursively enumerable sets of integer sequences which in some sense ``code'' the defining relations of groups. 

The algorithm we suggested in \cite{The Higman operations and  embeddings} allows  to list certain wide classes of groups for which Higman's famous embedding construction can also be constructive and explicit. 

\medskip
As a first step, a recursive group 
$G = \langle  A \mathrel{|} R \rangle = \langle a_1,a_2,\ldots \mathrel{|} r_1, r_2,\ldots \rangle $
with recursively enume\-rable relations $r_1, r_2,\ldots$
%($r_k=r_k(a_{i_1},\ldots,a_{i_k})$)
can be constructively embedded into a 
$2$-generator group
$T=\langle b,c 
\mathrel{|} r'_1, r'_2,\ldots 
\rangle$ 
where the relations $r'_1=r'_1(b,c),\;r'_2=r'_2(b,c),\ldots$ are certain words on just two letters $b,c$,  and they also are recursively enumerable (see  \cite{HigmanNeumannNeumann}
for the original embedding theorem and 
 \cite{Embeddings using universal words} for a method of embedding  that preserves the recursive enumeration).
Then for each $r'_n$,\; $n\!=\!1,2,\ldots$, a unique sequence $f_n$ of integers is compiled so that the set $\{r'_1, r'_2,\ldots \}$ of relations
is ``coded'' by means of the set $\mathcal B=\big\{f_1, f_2,\ldots \big\}$ of such sequences. 
Namely, for a relation: 
$$
r'_i  =r'_i(b,c) = b^{n_0}c^{n_1} \cdots b^{n_{2m}}c^{n_{2m+1}}
$$ 
for some $m=m(i)$, and $n_1,\ldots,n_{2m} \neq 0$ (the cases $n_0=0$, or $n_{2m+1}$ are \textit{not} ruled out) we output:
$$
f_i =(n_0, n_1, \ldots ,n_{2m+1}).
$$
Say, for the commutator word: 
$$r_1(b,c)=[b,c]=b^{-1}c^{-1}bc$$ 
we have the sequence: 
$$
f_1=(-1, -1, 1, 1).
$$

Since the transaction from relations set $R$ to sequences set $\mathcal B$ is done via a just few constructive steps, the set $\mathcal B$ also is recursively enumerable.

\medskip
Further for each $f_i\in \mathcal B$ Higman sets some special elements $b_{f_i}$ and $a_{f_i}$ in the free group $F_3 = \langle a,b,c\rangle$ of rank $3$, and using them defines certain respective subgroup $A_{\mathcal B}$ in $F_3$.

Then the so-called \textit{benign} subgroups are defined. One of the key results of \cite{Higman Subgroups of fP groups} is that \textit{$\mathcal B$ is recursively enumerable if and only if  $A_{\mathcal B}$ is benign in $F_3$}.

Finally ``the Higman Rope Trick'' in \cite{Higman Subgroups of fP groups} uses this benign subgroup $A_{\mathcal B}$ to embed $T$, and thus also the initial group $G$, into a finitely presented group.

\section{The explicit constructions}
\noindent
Our note \cite{The Higman operations and  embeddings} and this talk mainly concern the \textit{explicit construction of the recursive sets $\mathcal B$}, so we will discuss them in the sequel. 
The hard part of \cite{Higman Subgroups of fP groups}
is to show that $\mathcal B$ is recursively enumerable if and only if $\mathcal B$ can be constructed by some chain of special operators suggested by Higman:
\begin{equation}
\label{EQ Higman operations}
\iota,\; 
\upsilon,\; 
\rho,\; 
\sigma,\; 
\tau,\; 
\theta,\; 
\zeta,\; 
\pi,\; 
\omega_m.
\end{equation}
And parallel to application of those operations a respective \textit{benign subgroup} is being constructed in the free group $F_3$.

\medskip
In~\cite{Higman Subgroups of fP groups} Higman just relies on \textit{theoretical possibility} for construction of $\mathcal B$ via special operations \eqref{EQ Higman operations}, without any \textit{examples} of such construction for certain particular recursive groups. 
Is it worth noting that after Higman's result there was no attempt to explicitly find constructions of $\mathcal B$ by Higman's operations  for particular groups.
At least, we haven not seen them in the literature.

\medskip
We noticed that such a construction may be a doable task for some classes of groups for which the set $\mathcal B$ obeys certain simple ``combinatorial'' rules.
To explain what we understand under ``combinatorial'' rule let us bring an oversimplified example. Assume $\mathcal B$ consists of sequences $f_i\! =\!(n_0, n_1, \ldots ,n_{2m+1})$ where, say:\\
$n_0$ accepts any value grater then a fixed integer,  \\
$n_1$ accepts any value between two fixed integers, \\
$n_2$ accepts any value less than a fixed integer, \\
$n_3$ accepts arbitrary integer values, \\
$n_4$ and $n_5$ can accept any values but they are opposites of each other, etc. 
(an example will be given for $\mathbb Q$ below).

We noticed that this pattern does occur for many classes of groups, such as, the free abelian, metabelian, soluble, nilpotent groups, the additive group of  rational numbers $\mathbb Q$, the quasicyclic group $\mathbb C_{p^\infty}$, 
divisible abelian groups, etc.

Say, for the group $\mathbb Q$ by Example 3.5 in \cite{The Higman operations and  embeddings} the respective set $\mathcal B$ consists of all tuples $f$ of type:
$$
f_k=
\big(1,-k,-1, -k , -1,\; k,\;  1,\;  k, \; 1,\;  1, -1,\;  1\!\!-\!k, -1,-1,\;  1,\;  k\!-\!\!1,\;  1,\;  k\!-\!\!1,\;  -1\big).
$$
for $k=2,3,\ldots$, and so $\mathbb Q$ certainly is one of the groups which from the perspective of Higman embeddings has a not very complicated set of  relations.

Then we suggest \textit{an algorithm} with some generic tools that allow to explicitly construct $\mathcal B$ by Higman's operations for groups from the listed classes, and more.

\medskip
For this algorithm in addition to the initial Higman operations we suggest a few extra \textit{auxiliary operations} which make our work with the Higman operations more convenient and intuitive, see Subsection 2.4 in \cite{The Higman operations and  embeddings}.

\medskip
We intentionally brought our example above for the rational group $\mathbb Q$ because this topic is related to Problem 14.10 (a) posted by Bridson and de la Harpe in Kourovka Notebook as a ``well-known problem''.

Namely they asked to explicitly embed the group $\mathbb Q$ of rational numbers into a finitely presented group. 

\medskip
Our algorithm in \cite{The Higman operations and  embeddings} gives the method how such an explicit embedding can be build, without displaying any explicit embedding though.
An explicit embedding was recently found in in \cite{Belk Hyde Matucci} by Belk, Hyde and Matucci.  \cite{Belk Hyde Matucci} also stresses that \cite{The Higman operations and  embeddings} \textit{``has described how to explicitly carry out Higman's construction for $\mathbb Q$ as well as many other groups of interest''.}

\medskip
In \cite{The Higman operations and  embeddings} we intentionally try to stay as close to Higman's technique with free constructions of groups as possible, and also to the methods of  \cite{A modified proof for Higman}. Some of the proof details could be shortened by using tricks with wreath products of groups from  \cite{Subnormal embedding theorems},
 \cite{metabelian},
 \cite{On abelian subgroups of},
\cite{The Criterion of Shmel'kin},  
 \cite{Subvariety structures}, 
 \cite{Mikaelian 2002}, \cite{Mikaelian 2003}, \cite{Mikaelian 2003b}, \cite{Mikaelian 2004}, \cite{Mikaelian 2005}, \cite{Mikaelian 2007} or \cite{Mikaelian 2017},
however, we refrain from mixing free constructions technique with wreath product methods.

\end{document}